\documentclass[MSNbibl,number,citesort,dvips]{arxbj}
\usepackage{mathbh}

\aid{0}
\volume{18}
\issue{3}
\pubyear{2012}
\firstpage{803}
\lastpage{822}
\doi{10.3150/11-BEJ358}

\makeatletter
\newcommand{\U}{\mathbb{U}}
\newcommand{\V}{\mathbf{V}}
\newcommand{\N}{\mathbb{N}}
\newcommand{\M}{\mathbb{M}}
\newcommand{\R}{\mathbb{R}}
\newcommand{\C}{\mathbb{C}}
\newcommand{\F}{\mathbb{F}}
\newcommand{\D}{\mathbb{D}}
\newcommand{\I}{\mathbb{I}}
\newcommand{\oR}{\overline{\mathbb{R}}}
\newcommand{\pr}{\mathbb{P}}
\newcommand{\ex}{\mathbb{E}}
\newcommand{\covi}{\operatorname{\mathbb{C}ov}}
\newcommand{\eins}{\mathbh{1}}

\newcommand{\loc}{{\mathrm{loc}}}
\newcommand{\dis}{{\mathrm{d}}}

\newcommand{\BV}{\mathbb{BV}}

\newcommand{\cal}{\mathcal}
\newproclaim{assumption}{Assumption}[section]
\newremark{remarknorm}[assumption]{Remark}
\newtheorem{theorem}[assumption]{Theorem}
\newremark{examplenorm}[assumption]{Example}
\newtheorem{lemma}[assumption]{Lemma}
\makeatother

\begin{document}
\begin{frontmatter}

\title{Deriving the asymptotic distribution of\\ U- and V-statistics of
dependent data using weighted empirical processes}
\runtitle{Asymptotic distribution of U- and V-statistics}

\begin{aug}
\author[1]{\fnms{Eric} \snm{Beutner}\thanksref{1}\ead[label=e1]{e.beutner@maastrichtuniversity.nl}} \and
\author[2]{\fnms{Henryk} \snm{Z\"ahle}\corref{}\thanksref{2}\ead[label=e2]{zaehle@math.uni-sb.de}}
\runauthor{E.~Beutner and H.~Z\"ahle}
\address[1]{Department of Quantitative Economics, Maastricht
University, P.O.~Box 616, NL--6200 MD Maastricht, The Netherlands.
\printead{e1}}
\address[2]{Department of Mathematics, Saarland University, Postfach
151150, D--66041 Saarbr\"ucken, Germany.
\printead{e2}}
\end{aug}

\received{\smonth{7} \syear{2010}}
\revised{\smonth{11} \syear{2010}}

%
\begin{abstract}
It is commonly acknowledged that V-functionals with an unbounded kernel
are not Hadamard differentiable and that therefore the asymptotic
distribution of U- and V-statistics with an unbounded kernel cannot be
derived by the Functional Delta Method (FDM). However, in this article
we show that V-functionals are quasi-Hadamard differentiable and that
therefore a modified version of the FDM (introduced recently in
(\textit{J.~Multivariate Anal.} \textbf{101} (2010) 2452--2463)) can
be applied to this problem. The modified FDM
requires weak convergence of a weighted version of the underlying
empirical process. The latter is not problematic since there exist
several results on weighted empirical processes in the literature; see,
for example, (\textit{J.~Econometrics}
\textbf{130} (2006)
307--335,
\textit{Ann. Probab.}
\textbf{24} (1996)
2098--2127,
\textit{Empirical Processes with Applications to Statistics} (1986)
Wiley,
\textit{Statist. Sinica}
\textbf{18} (2008)
313--333). 
The modified FDM approach has the advantage that it is very flexible
w.r.t. both the underlying data and the estimator of the unknown
distribution function. Both will be demonstrated by various examples.
In particular, we will show that our FDM approach covers mainly all the
results known in literature for the asymptotic distribution of U- and
V-statistics based on dependent data -- and our assumptions are by
tendency even weaker. Moreover, using our FDM approach we extend these
results to dependence concepts that are not covered by the existing literature.
\end{abstract}

%

\end{frontmatter}

\section{Introduction}\label{Introduction}

For a distribution function (d.f.) $F$ on the real line, we consider the
characteristic
\begin{equation}\label{ucharacteristicdef}
U(F):=\int\int g(x_1,x_2)\,\mathrm{d}F(x_1)\,\mathrm{d}F(x_2)
\end{equation}
with $g\dvtx \R^2\to\R$ some measurable function, provided the double
integral exists. A systematic theory for the nonparametric estimation
of $U(F)$ was initiated in \cite{Hoeffding1948} and \cite
{vonMises1947}. A natural estimator for $U(F)$ is given by
\begin{equation}\label{ustatisticdef}
U(F_n):=\int\int g(x_1,x_2)\,\mathrm{d}F_n(x_1)\,\mathrm{d}F_n(x_2),
\end{equation}
where $F_n$ denotes some estimate of $F$ based on the first $n$
observations of a sequence $X_1,X_2,\ldots$ of random variables (on
some probability space $(\Omega,{\cal F},\pr)$) being identically
distributed according to $F$. Sometimes $U(F_n)$ is called
von-Mises-statistic (or simply V-statistic) with kernel $g$. If $F_n$
is the empirical d.f. $\hat F_n:=\frac{1}{n}\sum_{i=1}^n\eins
_{[X_i,\infty)}$ of $X_1,\ldots,X_n$, then we obtain
\begin{equation}\label{defUstat-sim}
U(\hat F_n)=\frac{1}{n^2}\sum_{i=1}^n\sum_{j=1}^ng(X_i,X_j),
\end{equation}
and we note that $U(\hat F_n)$ is closely related to the U-statistic
\begin{equation}\label{defUstat}
U_n:=\frac{1}{n(n-1)}\sum_{i=1}^n\sum_{j=1:j\not=i}^ng(X_i,X_j).
\end{equation}
If $X_1,\ldots,X_n$ are i.i.d., then $U_n$ is an unbiased estimator
whereas $U(\hat F_n)$ is generally not so. However, $U_n$ and $U(\hat
F_n)$ typically share the same asymptotic properties; cf. Remark~\ref
{mainthmholdsalsoforUS} below. Also notice that, in the nonparametric
setting, $U_n$ is the minimum variance unbiased estimator of $U(F)=\ex
[g(X_1,X_2)]$ whenever $X_1,\ldots,X_n$ are i.i.d. For background on
U-statistics see, for instance, \cite
{Dehling2006,Denker1985,Hoeffding1948,Lee1990,Sen1981,Sen1992,Serfling1980}.

We note that several features of a d.f. $F$ can be expressed as in~(\ref
{ucharacteristicdef}), for instance, the variance of $F$, or Gini's
mean difference of two independent random variables with d.f.~$F$; for
details, see Section~\ref{examples}.

Our objective is the asymptotic distribution of $U(F_n)$, that is, the
weak limit of the empirical error $\sqrt{n}(U(F_n)-U(F))$. In the
existing literature, the starting point for the derivation of the
asymptotic distribution of U-statistics $U_n$ is usually the Hoeffding
decomposition \cite{Hoeffding1948} of $U_n$. Using this
decomposition, asymptotic normality of $U_n$ was shown in \cite
{Hoeffding1948} for i.i.d. sequences, in \cite{Sen1972} for *-mixing
stationary sequences, in \cite{DenkerKeller1986,Yoshihara1976} for
$\beta$-mixing stationary sequences, in \cite{DewanRao2002} for
associated random variables, and recently in \cite{DehlingWendler2010}
for $\alpha$-mixing stationary sequences (recall from \cite{Bradley2005}, page 109: i.i.d. $\Rightarrow$ $*$-mixing $\Rightarrow$
$\beta$-mixing $\Rightarrow$ $\alpha$-mixing). Another approach is
based on the orthogonal expansion of the kernel $g$; see, for example,
\cite{DewanRao2001} and the references therein.

In this article, we derive the asymptotic distribution of U- and
V-statistics by means of a Functional Delta Method (FDM). The use of an
FDM is known to be beneficial for the following reason. Provided the
functional $U$ can be shown to be Hadamard differentiable at $F$, it is
basically enough to derive the asymptotic distribution of $F_n$ to
obtain the asymptotic distribution of $U(F_n)$. Therefore, this method
is especially useful for deriving the asymptotic distribution of the
estimator $U(\hat F_n)$ based on dependent data, because~-- given the
Hadamard differentiability -- one ``only'' has to derive the asymptotic
distribution of $\hat F_n$ based on data subjected to a certain
dependence structure. There are already several respective results on
the asymptotic distribution of $\hat F_n$ based on dependent data in
the literature (e.g., \cite{ChenFan2006,ShaoYu1996,Wu2008}), and new
respective results (combined with
the assumed Hadamard differentiability) would immediately yield also
the asymptotic distribution of~$U(\hat F_n)$.

However, one has to be careful with the application of an FDM to our
problem. The classical FDM in the sense of
\cite{Fernholz1983,Gill1989,Reeds1976} (see also \cite
{vanderVaart1998,vanderVaartWellner1996}) cannot be applied to many
interesting statistical functionals depending on the tails of the
underlying distribution, because the method typically relies on
Hadamard differentiability w.r.t. the uniform sup-norm. For
instance, as pointed out in \cite{vanderVaart1998} and \cite{Sen1996},
whenever $F$ has an unbounded support Hadamard differentiability
w.r.t. the uniform sup-norm can be shown neither for an L-statistic
with a weight function having one of the endpoints (or
both endpoints) of the closed interval $[0,1]$ in its support nor for a
U-statistic
with unbounded kernel. However, in \cite{BeutnerZaehle2010} a
modified version of the FDM was introduced which is suitable also
for nonuniform sup-norms (imposed on the tangential space only), and
it was in particular shown that this modified version can also be
applied to L-statistics with a weight function having one of the endpoints
(or both endpoints) of the closed interval
$[0,1]$ in its support. In contrast to the classical FDM, our FDM is
based on the
notion of quasi-Hadamard differentiability and requires weak
convergence of the empirical process $\sqrt{n}(\hat F_n-F)$ w.r.t. a
nonuniform sup-norm, that is, in other words, weak convergence of a
weighted version of the empirical process. Fortunately, the latter
is not problematic, because there are many results on the weak
convergence of weighted empirical processes in the literature; see
\cite{ShorackWellner1986} for i.i.d. data, and
\cite{ChenFan2006,ShaoYu1996,Wu2008} for dependent data.

In the present article, we demonstrate that the modified version of
the FDM can be applied to derive the limiting distribution for U- and
V-statistics with an unbounded kernel $g$. For simplicity of
notation, we restrict the derivations to kernels of degree 2.
However, in Remark~\ref{remarkond3}, we clarify how the results can
be extended to kernels of degree $d \geq3$. Using our FDM approach,
we will be able to a great extent to recover the results mentioned
above (the conditions imposed by our approach will turn out to be
weaker by tendency) and to extend them to other concepts of
dependence; cf. Section~\ref{examples-subs2}. The FDM approach
will also turn out to be useful when the empirical d.f. is replaced by
a different estimate of $F$, for instance by a smoothed version of
the empirical d.f.; cf.
Example~\ref{examplesmoothededfiid}.\looseness=1

The remainder of this article is organized as follows. In Section
\ref{sectionmainresult}, we state the conditions under which the
asymptotic distribution of U- and V-statistics can be derived by the
modified version of the FDM and present our main result. The
conditions imposed can be divided into two parts: on the one hand
conditions on the kernel $g$ and the d.f.~$F$, and on the other hand
conditions on an empirical process. In Section~\ref{examples}, we
give several examples for both, that is, for kernels $g$ and d.f. $F$ as
well as empirical processes fulfilling the conditions imposed.
In the Appendix~\ref{Jordandecomposition}, we recall the Jordan
decomposition of functions of locally bounded variation, which will
be beneficial for our applications in Section~\ref{examples}.
Finally, in the Appendix~\ref{appendixintegration} we give an
integration-by-parts formula and a sort of weighted Helly-Bray
theorem. Both results are needed in Section~\ref{sectionquasihadamard}
to show quasi-Hadamard differentiability of V-functionals.


\section{Main result}\label{sectionmainresult}

Our main result is Theorem~\ref{maintheorem} below, which provides a
CLT for the V-statistic $U(F_n)$ subject to Assumption \ref
{basicassumptions}. Let $\D_\lambda$ be the space of all  c\`adl\`ag
functions $\psi$ on $\oR$ with $\|\psi\|_\lambda<\infty$, where $\|
\psi\|_\lambda:=\|\psi\phi_\lambda\|_\infty$ refers to the
nonuniform sup-norm based on the weight function
$\phi_\lambda(x):=(1+|x|)^\lambda$, for $\lambda\in\R$ fixed. As
usual, we let $0 \cdot\infty:=0$. If $\lambda\ge0$, then we equip
$\D_\lambda$ with the $\sigma$-algebra ${\cal D}_\lambda:={\cal
D}\cap\D_\lambda$ to make it a measurable space, where~${\cal D}$ is
the $\sigma$-algebra generated by the usual coordinate projections
$\pi_x\dvtx \D\to\R$, $x\in\oR$, with~$\D$ the space of all bounded
 c\`adl\`ag  functions on $\oR$. Further, let $\BV_{\loc}$ be the space
of all functions $\psi\dvtx \oR\to\oR$ being real-valued and of local
bounded variation on $\R$. For $\psi\in\BV_{\loc}$, we denote by
$\mathrm{d}\psi^+$ and $\mathrm{d}\psi^-$ the unique positive Radon measures induced by
the Jordan decomposition of $\psi$ (for details, see the Appendix \ref
{Jordandecomposition}), and we set $|\mathrm{d}\psi|:=\mathrm{d}\psi^++\mathrm{d}\psi^-$.
Finally, we will interpret integrals as being over the open interval
$(-\infty,\infty)$, that is, $\int=\int_{(-\infty,\infty)}$.

\begin{assumption}\label{basicassumptions}
We assume that for some $\lambda>\lambda'\ge0$ the following
assertions hold:
\begin{enumerate}[(b)]
\item[(a)] For every $x_2\in\R$ fixed, the function $g_{x_2}(\cdot
):=g( \cdot ,x_2)$ lies in $\BV_\loc\cap\D_{-\lambda'}$.
Moreover, the function $x_2\mapsto\int\phi_{-\lambda}(x_1)
|\mathrm{d}g_{x_2}|(x_1)$ is measurable and finite w.r.t. $\|\cdot\|_{-\lambda'}$.
\item[(b)] The functions $g_{1,F}(\cdot):=\int g( \cdot
,x_2)\,\mathrm{d}F(x_2)$ and $g_{2,F}(\cdot):=\int g(x_1, \cdot )\,\mathrm{d}F(x_1)$ lie
in $\BV_{\loc}\cap \mathbb{D}$, and $\int\phi_{-\lambda}(x) |\mathrm{d}g_{i,F}|(x)<\infty
$ for $i=1,2$. Moreover, the functions $\overline{g_{1,F}}(\cdot
):=\int|g( \cdot , x_2)|\,\mathrm{d}F(x_2)$ and $\overline{g_{2,F}}(\cdot
):=\int|g(x_1, \cdot )|\,\mathrm{d}F(x_1)$ lie in $\D_{-\lambda'}$.
\item[(c)] $F$ is continuous, the double integral in (\ref
{ucharacteristicdef}) exists, and $\int\phi_{\lambda
'}(x)\,\mathrm{d}F(x)<\infty$. 
\item[(d)] $F_n:\Omega\to\D$ is $({\cal F},{\cal D})$-measurable,
and every realization of $F_n$ is nonnegative and nondecreasing, has
variation bounded by 1, the double integral in (\ref{ustatisticdef})
exists and $\int\phi_{\lambda'}(x)\,\mathrm{d}F_n(x)<\infty$, for every $n\in
\N$.
\item[(e)] The process $\sqrt{n}(F_n-F)$ is a random element of $(\D
_\lambda,{\cal D}_\lambda)$ for all $n\in\N$, and there is some
random element $B^\circ$ of $(\D_\lambda,{\cal D}_\lambda)$ with
continuous samples such that
\begin{equation}\label{conditiononBcirc}
\sqrt{n}(F_n-F)\stackrel{d}{\to}B^\circ  \qquad \mbox{in }(\D
_\lambda,{\cal D}_\lambda,\|\cdot\|_\lambda).
\end{equation}
\end{enumerate}
\end{assumption}

The assumptions (a) and (b) will allow us to prove quasi-Hadamard
differentiability of the functional $U$ (defined in (\ref
{ucharacteristicdef})) at $F$; see Section~\ref{sectionquasihadamard}.
At first glance, they seem to be awkward but in an application their
verification is often straightforward, see Section \ref
{examples-subs1}. To understand the meaning of conditions (a) and (b),
let us suppose that we want to derive the asymptotic distribution of U-
and V-statistics by means of the classical FDM in the sense of \cite
{Fernholz1983,Gill1989,Reeds1976}. Then we would have to prove Hadamard
differentiability of the functional $U$ given by (\ref
{ucharacteristicdef}) at~$F$. If $F$ has an unbounded support this
could be done by imposing Assumptions~\ref{basicassumptions}(a) and~(b) with $\lambda'=0$, that is, with the uniform sup-norm. Thus, as
pointed out in the Introduction, an application of the classical FDM
for the derivation of the asymptotic distribution of U- and
V-statistics would, inter alia, require a uniformly bounded kernel $g$
(cf. \cite{Sen1996}). On the other hand, the modified FDM only
requires that this boundedness holds w.r.t. the weaker nonuniform
sup-norm $\|\cdot\|_{-\lambda'}$ for some $\lambda' \geq0$.
\begin{remarknorm}\label{basicassumptionsremark}
Notice that
\begin{enumerate}[(b)$'$]
\item[(a)$'$] Assumption~\ref{basicassumptions}(a) could,
alternatively, be imposed on $g_{x_1}$ defined similar as~$g_{x_2}$.
Further notice that the second requirement in Assumption \ref
{basicassumptions}(a) is rather weak. Indeed: In the examples to be
given in Section~\ref{examples-subs1} the function $x_2\mapsto\int
\phi_{-\lambda}(x_1) |\mathrm{d}g_{x_2}|(x_1)$ even lies in~$\D$.
\item[(b)$'$] The last part of Assumption~\ref{basicassumptions}(b)
implies $g_{1,F},g_{2,F}\in\D_{-\lambda'}$.
\item[(c)$'$] Continuity of $F$ is required for the application of the
modified FDM.
\item[(d)$'$] Assumption~\ref{basicassumptions}(d) is always
fulfilled if $F_n$ is the empirical d.f. $\hat F_n$.
\item[(e)$'$] Assumption~\ref{basicassumptions}(e) does not require
that $F$ lies in $\D_\lambda$ or that $F_n$ is a random element of
$(\D_\lambda,{\cal D}_\lambda)$. These conditions would actually
fail to hold.
\end{enumerate}
\end{remarknorm}

\begin{theorem}\label{maintheorem}
Under Assumption~\ref{basicassumptions}, we have
\begin{equation}\label{convofempprocess}
\sqrt{n}\bigl(U(F_n)-U(F)\bigr)\stackrel{d}{\longrightarrow}\dot U(B^\circ
)  \qquad \mbox{in }(\R,{\cal B}(\R),|\cdot|)
\end{equation}
with
\begin{equation}\label{defoflimitofempprocess}
\dot U_F(B^\circ):=-\int B^\circ(x)\,\mathrm{d}g_{1,F}(x)-\int B^\circ(x)\,\mathrm{d}g_{2,F}(x).
\end{equation}
\end{theorem}

\begin{pf}
First of all, notice that the integrals in (\ref
{defoflimitofempprocess}) exist by Assumptions~\ref{basicassumptions}(b) and (e). Now, let $\BV_{1,\dis}$ be the space of all  c\`adl\`ag
functions in $\BV_{\loc}$ with variation bounded by~$1$, and $\U$ be
the class of all nonnegative and nondecreasing functions $f\in\BV
_{1,\dis}$ for which the integral on
the right-hand side of equation (\ref{defofU}) below and the integral $\int
\phi_{\lambda'}(x)\,\mathrm{d}f(x)$ exist. We define a functional $U\dvtx \U\to\R$
by setting
\begin{equation}\label{defofU}
U(f):=\int\int g(x_1,x_2)\,\mathrm{d}f(x_1)\,\mathrm{d}f(x_2), \qquad   f\in\U,
\end{equation}
so that $U(F)$ and $U(F_n)$ defined in (\ref
{ucharacteristicdef})--(\ref{ustatisticdef}) can be written as $U(f)$
with $f:=F$
and $f_n:=F_n$, respectively. We are going to apply an FDM to the
functional $U$. The version of the FDM we need for our purposes is
given in \cite{BeutnerZaehle2010}, Theorem 4.1. It is based on the
notion of quasi-Hadamard differentiability which is also introduced in
\cite{BeutnerZaehle2010}, Definition 2.1.

Let $\C_\lambda$ be the space of all continuous functions in $\D
_\lambda$, and notice that $\C_\lambda$ is separable w.r.t. $\|
\cdot\|_\lambda$. For every $f$ in $U$'s domain $\U$ we define a
functional $\dot U_f\dvtx \C_\lambda\to\R$ by setting
\begin{equation}\label{DefHadamard-functional}
\dot U_f(v) := -\int v(x)\,\mathrm{d}g_{1,f}(x)-\int v(x)\,\mathrm{d}g_{2,f}(x), \qquad   v\in
\C_\lambda,
\end{equation}
where $g_{i,f}$ is defined analogously to $g_{i,F}$ (cf. Assumption
\ref{basicassumptions}(b)). Lemma~\ref{HDofU} below shows that,
subject to Assumption~\ref{basicassumptions}(a)--(c), the functional
$U$ is quasi-Hadamard differentiable at $f:=F$ tangentially to $\C
_\lambda\langle\D_\lambda\rangle$ with quasi-Hadamard derivative
$\dot U_F$. Thus, assumption (iv) of Theorem~4.1 in \cite
{BeutnerZaehle2010} (with $f=U$, $V_f=\U$, $(\V',\| \cdot \|_{\V
'})=(\R,| \cdot |)$, $(\V_0,\| \cdot \|_{\V_0})=(\D_\lambda
,\| \cdot \|_{\lambda})$, $\C_0=\C_\lambda$, $\theta=F$ and
$T_n=F_n$) is fulfilled. Therefore, the statement of Theorem~\ref
{maintheorem} would follow from the FDM given in Theorem~4.1 in \cite
{BeutnerZaehle2010} if we could verify that\vadjust{\goodbreak} also the conditions
(i)--(iii) of this theorem are satisfied. Conditions (i) and (ii) are
satisfied by Assumption~\ref{basicassumptions}(d) and (e),
respectively. It thus remains to verify (iii), that is, that the
mapping $\widetilde\omega\mapsto U(W(\widetilde\omega)+F)$ is
$(\widetilde{\cal F},{\cal B}(\R))$-measurable whenever $W$ is a
measurable mapping from some measurable space $(\widetilde\Omega
,\widetilde{\cal F})$ to $(\D_\lambda,{\cal D}_\lambda)$ such that
$W(\widetilde\omega)+F\in\U$ for all $\widetilde\omega\in
\widetilde\Omega$. Since $W$ is $(\widetilde{\cal F},{\cal
D}_\lambda)$-measurable and ${\cal D}_\lambda$ is the projection
$\sigma$-field, we obtain in particular $(\widetilde{\cal F},{\cal
B}(\R))$-measurability of $\widetilde\omega\mapsto W(x,\widetilde
\omega)$ for every $x\in\oR$. Along with the representation (\ref
{ucharacteristicdef}), this yields $(\widetilde{\cal F},{\cal B}(\R
))$-measurability of $\widetilde\omega\mapsto
U(W(\widetilde\omega)+F)$.\vspace*{-3pt}
\end{pf}

We emphasize that Theorem~\ref{maintheorem} is quite a flexible
tool to derive the asymptotic distribution of the plug-in estimate
$U(F_n)$. In fact: Apart from checking the technical Assumptions
\ref{basicassumptions}(a)--(d), it is enough to establish the CLT
(\ref{conditiononBcirc}) for $F_n$ in order to obtain the CLT
(\ref{convofempprocess}) for $U(F_n)$. Section~\ref{examples}
below demonstrates this flexibility by various examples.

\begin{remarknorm}\label{maintheorem-remark}
If $B^\circ$ in Theorem~\ref{maintheorem} is a Gaussian process with
zero mean and measurable covariance function $\Gamma$ and if $\int
\int\Gamma(x,y)\,\mathrm{d}g_{i,F}(x)\,\mathrm{d}g_{j,F}(y)$ exists for every $i,j\in\{
1,2\}$, then the random variable $\dot U_F(B^\circ)$ defined in (\ref
{defoflimitofempprocess}) is normally distributed with mean $0$ and variance
\begin{equation}\label{asymptvar}
\sigma^2 := \sum_{i=1}^2\sum_{j=1}^2 \int\int\Gamma(x,y)\,\mathrm{d}g_{i,F}(x)\,\mathrm{d}g_{j,F}(y).
\end{equation}
\end{remarknorm}

\begin{remarknorm}\label{mainthmholdsalsoforUS}
If $\ex[|g(X_1,X_1)|]<\infty$ (in Examples~\ref{ginimeandifference}
and~\ref{examplevariance} below we even have $g(x,x)=0$ for all $x\in
\R$), then the particular V-statistic $U(\hat F_n)$ and the
U-statistic $U_n$ (defined in (\ref{defUstat-sim}) and (\ref
{defUstat}), resp.) have the same asymptotic distribution. To
see this, we first of all note that (for $n\ge2$)
\begin{eqnarray}\label{equivofuandv}
&&\sqrt{n}\bigl(U_n-U(F)\bigr)\nonumber\\
&& \quad  =  \sqrt{n}\bigl(U_n-U(\hat F_n)\bigr)+\sqrt{n}\bigl(U(\hat
F_n)-U(F)\bigr)\nonumber
\\[-8pt]
\\[-8pt]
&& \quad  =  \frac{\sqrt{n}}{n-1} U(\hat F_n)-\frac{\sqrt{n}}{n(n-1)}\sum
_{i=1}^ng(X_i,X_i)+\sqrt{n}\bigl(U(\hat F_n)-U(F)\bigr)\nonumber\\
&& \quad  =:  S_1(n)-S_2(n)+\sqrt{n}\bigl(U(\hat F_n)-U(F)\bigr).
\nonumber
\end{eqnarray}
As $\sqrt{n}(U(\hat F_n)-U(F))$ converges weakly to some nondegenerate
limit, we obtain by Slutzky's lemma that $S_1(n)=\frac{1}{n-1} \sqrt
{n}(U(\hat F_n)-U(F))+\frac{\sqrt{n}}{n-1} U(F)$ converges in
probability to zero. Further, by the Markov inequality we know that,
for every $\varepsilon>0$ fixed, $\pr[|S_2(n)|>\varepsilon]$ is
bounded above by $\frac{1}{\varepsilon}\ex[|S_2(n)|]$ which, in
turn, is bounded above by $\frac{\sqrt{n}}{n-1}\frac{1}{\varepsilon
}\ex[|g(X_1,X_1)|]$. So we also have that $S_2(n)$ converges in
probability to zero. Slutzky's lemma and (\ref{equivofuandv}) thus
imply that $\sqrt{n}(U_n-U(F))$ has indeed the same limit distribution
as $\sqrt{n}(U(\hat F_n)-U(F))$.
\end{remarknorm}

\begin{remarknorm}
The linear part of the Hoeffding decomposition of $U_n-U(F)$ (cf.
\cite{Serfling1980}, page~178) multiplied by $\sqrt{n}$ can be written
as $\sum_{i=1}^2\int g_{i,F}\,\mathrm{d}(\sqrt{n}(\hat F_n-F))$, for example, using\vadjust{\goodbreak}
the integration-by-parts formula (\ref{defofintegraludv-2}), as $-\sum
_{i=1}^2\int\sqrt{n}(\hat F_n-F)\,\mathrm{d}g_{i,F}$. Then, if we could show
that the degenerate part of $U_n$ converges in probability to zero
(which is nontrivial for dependent data), we could recover (\ref
{convofempprocess}) with $U_n$ in place of $U(F_n)$ by using~(\ref
{conditiononBcirc}) and the Continuous Mapping theorem.\vspace*{-3pt}
\end{remarknorm}


\section{Examples}\label{examples}

In this section, we give some examples for $g$, $F$ and $F_n$
satisfying Assumption~\ref{basicassumptions}. At first, in Section
\ref{examples-subs1}, we provide examples for $g$ (and $F$) satisfying
\mbox{Assumptions~\ref{basicassumptions}(a)--(b)}. Thereafter, in
Section~\ref{examples-subs2}, we will give examples for $F_n$ (and
$F$) satisfying Assumptions~\ref{basicassumptions}(d)--(e) for
various types of data. We assume throughout this section that
Assumption~\ref{basicassumptions}(c) is fulfilled because its
meaning is rather obvious and the conditions imposed by it are fairly
weak.\vspace*{-3pt}


\subsection{Examples for $g$}\label{examples-subs1}

In \cite{Aaronsonetal1996}, one can find a number of examples for
kernels $g$ for which $U(F)$ corresponds to a popular characteristic of
$F$. By means of two popular examples, we now illustrate how to verify
the Assumptions~\ref{basicassumptions}(a)--(b). It will be seen that
the verification of these assumptions is easy, though, at first glance,
it may seem cumbersome. We will use the notion of Jordan decomposition
$\psi=\psi(c)+\psi_c^+-\psi_c^-$ centered at some point $c\in\R$.
For the reader's convenience, we have recalled the essentials in the
Appendix~\ref{Jordandecomposition}.


\begin{examplenorm}[(Gini's mean difference)] \label{ginimeandifference}
If $g(x_1,x_2)=|x_1-x_2|$ and $F$ has a finite first moment, then
$U(F)$ equals Gini's mean difference $\ex[|X_1-X_2|]$ of two i.i.d.
random variables $X_1$ and $X_2$ on some probability space $(\Omega
,{\cal F},\pr)$ with d.f. $F$. Then the Assumptions \ref
{basicassumptions}(a)--(b) are fulfilled for $\lambda'=1$. Indeed:
We have $g_{x_2}(x_1)=(x_1-x_2)\eins_{(x_2,\infty
]}(x_1)-(x_1-x_2)\eins_{[-\infty,x_2]}(x_1)$, so that the first part
of Assumption~\ref{basicassumptions}(a) obviously holds. Further,
the Jordan decomposition (\ref{eqJordanwitha}) of $g_{x_2}$ centered
at $c=x_2$ reads as
$g_{x_2}(x_1)=0+{g_{x_2}}_{x_2}^+(x_1)-{g_{x_2}}_{x_2}^-(x_1)$, where
${g_{x_2}}_{x_2}^+(x_1)= (x_1-x_2) \eins_{(x_2,\infty]}(x_1)$ and
${g_{x_2}}_{x_2}^-(x_1)=(x_1-x_2) \eins_{[-\infty,x_2]}(x_1)$, and
so, in view of Lemma~\ref{independenceofc}, $\mathrm{d}g_{x_2}^+(x_1) = \eins
_{(x_2,\infty]}(x_1)\,\mathrm{d}x_1$ and $\mathrm{d}g_{x_2}^-(x_1) = \eins_{[-\infty
,x_2]}(x_1)\,\mathrm{d}x_1$. Now it can be seen easily that also the second part
of Assumption~\ref{basicassumptions}(a) holds; we omit the details.
Let us now turn to Assumption~\ref{basicassumptions}(b). We have
\begin{eqnarray*}
&& g_{1,F}(x_1)\\
&& \quad  =  \ex\bigl[X_2\eins_{(x_1,\infty]}(X_2)\bigr]-x_1\pr[X_2>x_1]+x_1\pr
[X_2\le x_1]-\ex\bigl[X_2\eins_{[-\infty,x_1]}(X_2)\bigr] \\
&& \quad  =  x_1\bigl(2F(x_1)-1\bigr)-\ex[X_2] +2 \ex\bigl[X_2\eins_{(x_1,\infty]}(X_2)\bigr]\\
&& \quad  =  K+x_1+2 \bigl(-x_1\bigl(1-F(x_1)\bigr)+ \ex\bigl[X_2\eins_{(x_1,\infty
]}(X_2)\bigr] \bigr)\\
&& \quad  =  K+x_1+2\int_{x_1}^{\infty}\bigl(1-F(x)\bigr)\,\mathrm{d}x
\end{eqnarray*}
with $K:=-\ex[X_2]$. The same representation holds for $g_{2,F}$. So
we obviously have $g_{i,F}=\overline{g_{i,F}}\in\D_{-1}\cap\BV
_\loc$ for $i=1,2$. Moreover,\vadjust{\goodbreak} we have $g_{i,F}'(x)=2F(x)-1$, and so
there is some constant $c\in\R$ such that $g_{i,F}$ is nonincreasing
on $(-\infty, c)$ and is nondecreasing on $(c,\infty)$, for $i=1,2$.
Since the density of $|\mathrm{d}g_{i,F}|$ on $(-\infty,c)$ and the density of
$|\mathrm{d}g_{i,F}|$ on $(c,\infty)$ are bounded, we also have $\int\phi
_{-\lambda}(x)|\mathrm{d}g_{i,F}|(x)<\infty$ for $i=1,2$ and every $\lambda
>1$. That is, all parts of Assumption~\ref{basicassumptions}(b) hold
true. Thus, Assumptions~\ref{basicassumptions}(a)--(b) hold~true.
\end{examplenorm}

If also Assumptions~\ref{basicassumptions}(d)--(e) hold true, then
we obtain from Theorem~\ref{maintheorem} for the kernel
$g(x_1,x_2)=|x_1-x_2|$ that $\dot U(B^\circ)= 2\int B^\circ
(x)(1-2F(x))\,\mathrm{d}x$, because $\mathrm{d}g_{1,F}(x)=\mathrm{d}g_{2,F}(x)=(2F(x)-1)\,\mathrm{d}x$.


\begin{examplenorm}[(Variance)]\label{examplevariance}
If $g(x_1,x_2)=\frac{1}{2}(x_1-x_2)^2$ and $F$ has a finite second
moment, then $U(F)$ equals the variance of $F$. In this case, the
Assumptions~\ref{basicassumptions}(a)--(b) are fulfilled for
$\lambda'=2$. The verification of this is even easier than the
elaborations in Example~\ref{ginimeandifference}. We note that this
time, we obtain $\mathrm{d}g_{x_2}^+(x_1)= (x_1-x_2) \eins_{(x_2,\infty
]}(x_1)\,\mathrm{d}x_1$ and $\mathrm{d}g_{x_2}^-(x_1) = (x_2-x_1) \eins_{[-\infty
,x_2]}(x_1)\,\mathrm{d}x_1$ as well as
$\mathrm{d}g_{i,F}^+(x_i) = (x_i-\ex[X_j]) \eins_{(\ex[X_j],\infty]}(x_i)\,\mathrm{d}x_i$ and $\mathrm{d}g_{i,F}^-(x_i) = (\ex[X_j]-x_i) \eins_{[-\infty,\ex
[X_j]]}(x_i)\,\mathrm{d}x_i$ for $i,j\in\{1,2\}$ with $i\not=j$.
\end{examplenorm}

If also Assumptions~\ref{basicassumptions}(d)--(e) hold true, then
we obtain from Theorem~\ref{maintheorem} for the kernel
$g(x_1,x_2)=\frac{1}{2}(x_1-x_2)^2$ that $\dot U(B^\circ)= 2\int
B^\circ(x)(\ex[X_1]-x)\,\mathrm{d}x$, because $\mathrm{d}g_{1,F}(x)=\mathrm{d}g_{2,F}(x)=(x-\ex
[X_1])\,\mathrm{d}x$.


\subsection{Examples for $F_n$}\label{examples-subs2}

Here we will give some examples for estimators $F_n$ for $F$ that
satisfy \mbox{Assumption~\ref{basicassumptions}(d)--(e)}. We first consider
the case of i.i.d. data.

\begin{examplenorm}[(Empirical d.f. of i.i.d. data)]\label{exampleedfiid}
Let $X_1,X_2,\ldots$ be a sequence of i.i.d. random variables with d.f.
$F$, and let $\lambda\geq0$. If $F$ has a finite $\gamma$-moment for
some $\gamma>2\lambda$, then Theorem 6.2.1 in \cite
{ShorackWellner1986} shows that for the empirical d.f. $\hat F_n$ of
$X_1,\ldots,X_n$,
\begin{equation}\label{exampleedfiid-CLT}
\sqrt{n}(\hat F_n-F)\stackrel{d}{\to} B_F^\circ  \qquad \mbox{(in
$(\D_\lambda,{\cal D}_\lambda,\|\cdot\|_\lambda)$)},
\end{equation}
where $B_F^\circ$ is an $F$-Brownian bridge, that is, a centered
Gaussian process with covariance function $\Gamma(x,y)=F(x\wedge
y)\overline{F}(x\vee y)$. Thus, if $\lambda>0$, if $F$ has a finite
$\gamma$-moment for some $\gamma>2\lambda$, and if $g$ is a kernel
satisfying Assumptions~\ref{basicassumptions}(a)--(b) for $F$ and
some $\lambda'\in[0,\lambda)$, then Theorem~\ref{maintheorem} shows
that the law of $\sqrt{n}(U(\hat F_n)-U(F))$ converges weakly to the
normal distribution with mean $0$ and variance given by (\ref
{asymptvar}) with $\Gamma(x,y)=F(x\wedge y)\overline{F}(x\vee y)$.
Alternatively, the result can be stated as follows: If $g$ is a
fixed kernel and $\mathbb{F}_{g, \lambda'}$ denotes the class of all
d.f. $F$ for which Assumptions~\ref{basicassumptions}(a)--(b) hold
with $\lambda' \geq0$, then $\sqrt{n}(U(\hat F_n)-U(F))$ converges
weakly to the above mentioned normal distribution for every
$F\in\mathbb{F}_{g, \lambda'}$ having a finite $\gamma$-moment for
some $\gamma>2 \lambda'$. Indeed: In this case, we can choose
$\lambda\in(\lambda',\gamma/2)$ in Assumption \ref
{basicassumptions}(e).
\end{examplenorm}


\begin{examplenorm}[(Smoothed empirical d.f. of i.i.d.
data)]\label{examplesmoothededfiid}
Suppose that in the setting of Example~\ref{exampleedfiid} the
empirical d.f. $\hat F_n$ is smoothed out by the heat kernel
$p_{\varepsilon_n}(\cdot)$ with bandwidth $\varepsilon_n\ge0$, that
is, that $\hat F_n$ is replaced by $P_{\varepsilon_n}\hat F_n$ with
$(P_\varepsilon)_{\varepsilon\ge0}$ the heat semigroup (i.e.,
$P_\varepsilon\psi:=\int_\R\psi(y)p_\varepsilon(\cdot-y)\,\mathrm{d}y$ for
$\varepsilon>0$, and $P_0:=\I$). Then, if $F$ is also Lipschitz
continuous and $\sqrt{n} \varepsilon_n^{(\gamma-\lambda)/(2\gamma
)}\to0$, the CLT (\ref{exampleedfiid-CLT}) (with $\hat F_n$ replaced
by $P_{\varepsilon_n}\hat F_n$) still holds (cf. Corollary A.2 in
\cite{BeutnerZaehle2010}), and therefore the weak limit of the law of
$\sqrt{n}(U(P_{\varepsilon_n}\hat F_n)-U(F))$ is still the normal
distribution with mean $0$ and variance given by (\ref{asymptvar})
with $\Gamma(x,y)=F(x\wedge y)\overline{F}(x\vee y)$. Of course, at
this point we have to ensure that under the imposed assumptions the
expression $U(P_{\varepsilon_n}\hat F_n)$ is well defined, that is,
that Assumption~\ref{basicassumptions}(d) is satisfied. Now, it can
be easily deduced from Lemma 3.2 in \cite{Zaehle2010} that in our
setting $P_{\varepsilon_n}\hat F_n$ lies in $\D_\lambda$. Thus, if
we assume that, for example, $\sup_{x_1,x_2\in\R}|g(x_1,x_2)|\phi
_{-\lambda'}(x_1)\phi_{-\lambda'}(x_2)<\infty$, Assumption \ref
{basicassumptions}(d) follows easily.
\end{examplenorm}


Let us now turn to the case of dependent data, which is our actual
objective. Throughout the examples presented below, we consider a
strictly stationary sequence $(X_i)=(X_i)_{i\ge1}$ of random
variables on some probability space $(\Omega,{\cal F},\pr)$ with
continuous d.f. $F$, and let as before $\hat F_n$ denote the
corresponding empirical d.f. at stage $n$. By strict stationarity, we
mean that the joint distribution of $X_{i+1},\ldots,X_{i+m}$ does
not depend on $i$ for every fixed positive integer $m$. We will
consider three popular dependency structures ($\alpha$-, $\beta$- and
$\rho$-mixing) in more detail in Examples~\ref{alphamixing},
\ref{betamixing}, and~\ref{rhomixing}, respectively. There, we
will also provide a comparison of the results obtained by the
approach considered here and the results obtained up to now. For the
definition of $\alpha$-, $\beta$- and $\rho$-mixing (and other)
mixing conditions and for examples of strictly stationary $\alpha$-,
$\beta$- and $\rho$-mixing sequences see, for example,
\cite{Bradley2005,Doukhan1994,Lindner2008}. As usual, the
corresponding mixing coefficients will be referred to as
$\alpha(n)$, $\beta(n)$ and $\rho(n)$, respectively. The application
of our method to other dependence concepts will be discussed in
Example~\ref{exampleWu+associated}. Notice that the condition of
$\alpha$-mixing is weaker than the condition of $\beta$-mixing
(absolute regularity) under which CLTs for U-statistics have been
established in \cite{DenkerKeller1986,Yoshihara1976}. A CLT for
strictly stationary $\alpha$-mixing (strongly mixing) sequences of
random variables has been given in \cite{DehlingWendler2010}.

\begin{examplenorm}[(Empirical d.f. of $\alpha$-mixing
data)]\label{alphamixing}
Let $(X_i)$ be $\alpha$-mixing with $\alpha(n)=\mathcal
{O}(n^{-\theta})$ for some $\theta>1+\sqrt{2}$, and let $\lambda\ge
0$. If $F$ has a finite $\gamma$-moment for some $\gamma>\frac
{2\theta\lambda}{\theta-1}$, then it can easily be deduced from
Theorem 2.2 in \cite{ShaoYu1996} that
\begin{equation}\label{example-dependent-1}
\sqrt{n}(\hat F_{n}-F)\stackrel{d}{\to}\widetilde B_F^\circ
 \qquad \mbox{(in $(\D_\lambda,{\cal D}_\lambda,\|\cdot\|_\lambda)$)}
\end{equation}
with $\widetilde B_F^\circ$ a continuous centered Gaussian process
with covariance function
\begin{eqnarray}
\label{gammafordependent-1}
\Gamma(s,t) & = &
F(s\wedge t)\bar F(s\vee t) \nonumber
\\[-8pt]
\\[-8pt]
& &{} +\sum_{k=2}^{\infty}\bigl [ \covi \bigl(\eins_{\{X_1 \leq s\}},
\eins_{\{X_k \leq t\}} \bigr)+ \covi \bigl(\eins_{\{X_1 \leq t\}},
\eins_{\{X_k \leq s\}} \bigr) \bigr]
\nonumber
\end{eqnarray}
(cf. Section 3.3 in \cite{BeutnerZaehle2010}). Thus, if $g$ is a
fixed kernel and $\mathbb{F}_{g, \lambda'}$ denotes the class of all
d.f. satisfying Assumptions~\ref{basicassumptions}(a)--(b) for some
$\lambda' \geq0$, then Theorem~\ref{maintheorem} shows that the
law of $\sqrt{n}(U(\hat F_n)-U(F))$ converges weakly to the normal
distribution with mean~$0$ and variance given by (\ref{asymptvar}),
with $\Gamma$ as in (\ref{gammafordependent-1}), for every d.f.
$F\in\F_{g, \lambda'}$ having a finite $\gamma$-moment for some
$\gamma> \frac{2\theta\lambda'}{\theta-1}$. Indeed: In this case
we can choose $\lambda\in(\lambda',\gamma(\theta-1)/(2\theta))$ in
Assumption~\ref{basicassumptions}(e).

To compare our result with that of Theorem 1.8 in
\cite{DehlingWendler2010}, we consider the kernel
$g(x_1,x_2)=\frac{1}{2}(x_1-x_2)^2$. For Theorem 1.8 in
\cite{DehlingWendler2010} to be applicable, we must assume that~$F$
has a finite $\gamma$-moment for some $\gamma>4$ (the same condition
is necessary to ensure that the approach considered here works). In
this case, both integrability conditions in Theorem 1.8 in
\cite{DehlingWendler2010} are fulfilled, and the condition on the
mixing coefficients reads as follows:
$\alpha(n)=\mathcal{O}(n^{-\theta})$ for some
$\theta>\frac{3}{2}+\frac{1}{2\gamma}+\frac{5}{\gamma-4}+\frac
{2}{\gamma(\gamma-4)}=\frac{3\gamma-1}{2\gamma-8}$.
On the other hand, if $F$ has a finite $\gamma$-moment for some
$\gamma>4$, in our setting we may choose $\lambda'=2$, and so
$\theta>\frac{\gamma}{\gamma-4}$ (and
$\lambda\in(2,\frac{\gamma(\theta-1)}{2\theta})$). Hence, our
condition on the mixing coefficients reads as follows:
$\alpha(n)=\mathcal{O}(n^{-\theta})$ for some
$\theta>\frac{\gamma}{\gamma-4}$. Notice that
$\frac{3\gamma-1}{2\gamma-8}>\frac{\gamma}{\gamma-4}$ holds for all
$\gamma>4$. Taking into account that in our setting, we must choose
$\theta> 1 + \sqrt{2}$ for the result of \cite{ShaoYu1996} to be
applicable we find that our result relies on a weaker assumption on the
mixing coefficients than Theorem 1.8 in \cite{DehlingWendler2010}
whenever $\frac{3\gamma-1}{2\gamma-8} > 1+\sqrt{2}$, that is,
$\gamma< \frac{7+8\sqrt{2}}{2\sqrt{2}-1}$.
\end{examplenorm}

\begin{examplenorm}[(Empirical d.f. of $\beta$-mixing data)]\label{betamixing}
Let $(X_i)$ be $\beta$-mixing with $\beta(n)=\mathrm{O}(n^{-\theta})$ for
some $\theta> \frac{\kappa}{\kappa-1}$ with $\kappa>1$, and let
$\lambda\ge0$. If $F$ has a finite $\gamma$-moment for some
$\gamma>2\lambda\kappa$, then it can easily be deduced from Lemma
4.1 in \cite{ChenFan2006} that the CLT (\ref{example-dependent-1})
still holds and that the covariance function is again given by
(\ref{gammafordependent-1}). Thus, if $g$ is a fixed kernel and
$\mathbb{F}_{g, \lambda'}$ denotes the class of all d.f. satisfying
Assumptions~\ref{basicassumptions}(a)--(b) for some $\lambda'
\geq0$, then Theorem~\ref{maintheorem} shows that the law of
$\sqrt{n}(U(\hat F_n)-U(F))$ converges weakly to the normal
distribution with mean $0$ and variance given by (\ref{asymptvar}),
with $\Gamma$ as in (\ref{gammafordependent-1}), for every d.f.
$F\in\F_{g, \lambda'}$ having a finite $\gamma$-moment for some
$\gamma>2 \lambda' \kappa$. Indeed: In this case we can choose
$\lambda\in(\lambda',\frac{\gamma}{2\kappa})$.

To compare our result with that of Theorem 3.1 in
\cite{Yoshihara1976} (see also Theorem 1.8 in
\cite{DehlingWendler2010}), we consider the kernel
$g(x_1,x_2)=\frac{1}{2}(x_1-x_2)^2$. For this theorem to be
applicable, we must again assume that $F$ has a finite
$\gamma$-moment for some $\gamma>4$ (the same condition is again
necessary to ensure that the approach considered here works). In
this case, both integrability conditions in Theorem 3.1 in
\cite{Yoshihara1976} (see also Theorem 1.8 in
\cite{DehlingWendler2010}) are fulfilled, and the condition on the
mixing coefficients reads as follows:
$\beta(n)=\mathcal{O}(n^{-\theta})$ for some
$\theta>\frac{\gamma}{\gamma-4}$. On the other hand, if $F$ has a
finite $\gamma$-moment for some $\gamma>4$, in our setting we may
choose $\lambda'=2$, and so $\kappa<\gamma/4$ (and
$\lambda\in(2,\frac{\gamma}{2\kappa})$). Hence, in view of
$\theta>\frac{\kappa}{\kappa-1}$, our condition on the mixing
coefficients reads as follows: $\beta(n)=\mathcal{O}(n^{-\theta})$
for some $\theta>\frac{\gamma}{\gamma-4}$. That is, both results
impose the same condition on the mixing coefficients.
\end{examplenorm}


\begin{examplenorm}[(Empirical d.f. of $\rho$-mixing data)]\label{rhomixing}
Let $(X_i)$ be $\rho$-mixing with $\sum_{n=1}^{\infty} \rho(2^n) <
\infty$, suppose
$
\sum_{k=2}^{\infty} |\covi(\eins_{\{X_1 \leq s\}}, \eins_{\{
X_k \leq t\}})+ \covi(\eins_{\{X_1 \leq t\}}, \eins_{\{X_k \leq s\}
}) |<\infty,
$
and let $\lambda\geq0$. If $F$ has a finite $\gamma$-moment for some
$\gamma>\lambda(2+\varepsilon)$ with $\varepsilon> 0$, then it can
easily be deduced from Theorem~2.3 in \cite{ShaoYu1996} that the CLT
(\ref{example-dependent-1}) still holds and that the covariance
function is again given by (\ref{gammafordependent-1}) (cf. Section
3.3 in \cite{BeutnerZaehle2010}). Hence, we again have in this case:
If $g$ is a fixed kernel and if we denote by $\mathbb{F}_{g, \lambda
'}$ the class of all d.f. for which Assumptions~\ref{basicassumptions}(a)--(b) hold for some $\lambda' \geq0$, then Theorem~\ref
{maintheorem} yields that the law of $\sqrt{n}(U(\hat F_n)-U(F))$
converges weakly to the normal distribution with mean $0$ and variance
given by (\ref{asymptvar}) with $\Gamma$ as in (\ref
{gammafordependent-1}) for every $F \in\mathbb{F}_{g, \lambda'}$
having a finite $\gamma$-moment for some $\gamma>\lambda
'(2+\varepsilon)$. Indeed: In this case, we can choose $\lambda\in
(\lambda',\gamma/(2+\varepsilon))$.
\end{examplenorm}

Up to our best knowledge, the asymptotic distribution of
U- and V-statistics of $\rho$-mixing data has not been studied explicitly
so far. Of course, every $\rho$-mixing sequence is also
$\alpha$-mixing (since $\alpha(n)\le\frac{1}{4}\rho(n)$; see
\cite{Bradley2005}, Inequality (1.12)), but the condition on the
mixing coefficients imposed in Example~\ref{rhomixing} is
considerably weaker than the condition on the mixing coefficients
imposed in Example~\ref{alphamixing}. Similar statements apply to
further dependence concepts, and one also obtains that further
dependence concepts are also covered by our approach.


\begin{examplenorm}[(Further examples)]\label{exampleWu+associated}
Recently, a new dependence structure for sequences of random
variables was introduced in \cite{Wu2008}. Thus, not surprising,
limit distributions for U- and V-statistics under this dependence
concepts have not been derived so far. Anyhow, in \cite{Wu2008} it
was also proved that, subject to certain conditions, the weighted
empirical process $\sqrt{n}(\hat F_n-F)\phi_\gamma$ converges weakly
to a
tight Gaussian process. Here $\hat F_n$ is the empirical d.f. based on a
sequence of random variables fulfilling this dependence condition.
From our Theorem~\ref{maintheorem} one can thus (along the lines of
Examples~\ref{alphamixing},~\ref{betamixing}, and~\ref{rhomixing})
derive the limit distribution of U- and V-statistics when
the data fulfills the dependence structure in \cite{Wu2008}. We omit
the details.

In \cite{DewanRao2002}, the limit distribution of U-statistics for
associated sequences was derived using the Hoeffding decomposition.
To prove asymptotic normality of U-statistics for stationary and
associated sequences, it was required there that the partial
derivatives of $g$ are uniformly bounded. This clearly excludes the
variance of a random variable. On the other hand, our approach also
covers the variance for the case of stationary and associated
sequences. Indeed: Let $(X_i)$ be a stationary, associated sequence
with $\covi(X_1,X_n)=\mathrm{O}(n^{-\nu-\varepsilon})$ for some $\nu\geq(3
+ \sqrt{33})\slash2$ and $\varepsilon> 0$. Then, we can deduce
from Theorem 2.4 in \cite{ShaoYu1996} that the CLT
(\ref{example-dependent-1}) still holds and the covariance function is
again given by (\ref{gammafordependent-1}) whenever $F$ has a finite
$\gamma$-moment for some $\gamma>\frac{2 \lambda\nu}{\nu-3}$
($\lambda\geq0$ fixed). Hence, we obtain from Theorem \ref
{maintheorem} (recall from Example~\ref{examplevariance} that
Assumptions~\ref{basicassumptions}(a)--(b) are fulfilled for the
variance with $\lambda'=2$) that the variance is included in our
method of proof whenever $F$ has a finite $\gamma$-moment for some
$\gamma>\frac{4 \nu}{\nu-3}$; in this case we can choose $\lambda
\in(2,\gamma(\nu-3)/(2\nu))$.
\end{examplenorm}


\section{Quasi-Hadamard differentiability of $U$}\label{sectionquasihadamard} 

This section is concerned with the quasi-Hadamard differentiability (in
the sense of Definition~2.1 in \cite{BeutnerZaehle2010}) of the
functional $U$ defined in (\ref{defofU}). Recall that quasi-Hadamard
differentiability is needed in the proof of Theorem~\ref{maintheorem}.
Recall also that $\BV_{1,\dis}$ is the space of all  c\`adl\`ag
functions in $\BV_{\loc}$ with variation bounded by $1$, and that $\U
$ is the class of all nonnegative and nondecreasing functions $f\in\BV
_{1,\dis}$ for which the integral on the right-hand side of equation (\ref
{defofU}) and the integral $\int\phi_{\lambda'}(x)\,\mathrm{d}f(x)$ exist.
Moreover, we let~$\BV_{\loc,\dis}$ be the space of all  c\`adl\`ag
functions in $\BV_{\loc}$.

\begin{lemma}\label{HDofU}
Under Assumptions~\ref{basicassumptions}\textup{(a)--(c)} (the continuity of
$F$ is actually superfluous at this point), the functional $U$ defined
in (\ref{defofU}) is quasi-Hadamard differentiable at $f:=F$
tangentially to $\C_\lambda\langle\D_\lambda\rangle$ with
quasi-Hadamard derivative given by $\dot U_f$ defined in (\ref
{DefHadamard-functional}) with $f:=F$.
\end{lemma}

\begin{pf}
To prove the claim, we have to show that
\begin{eqnarray}\label{defeqforHD}
\lim_{n\to\infty} \biggl| \dot U_f (v)-\frac
{U(f+h_nv_n)-U(f)}{h_n} \biggr|=0
\end{eqnarray}
holds for each triplet $(v,(v_n),(h_n))$ with $v\in\C_\lambda$,
$(v_n)\subset\D_\lambda$ satisfying $f+h_n v_n\in\U$ (for all
$n\in\N$) as well as $\|v_n -v\|_\lambda\to0$, and $(h_n)\subset\R
_0:=\R\setminus\{0\}$ satisfying $h_n\to0$. Let $f_n:=f+h_nv_n$. We
stress the fact that $f_n$ lies in $\U$ which is
a subset of $\BV_{1,\dis}$, and that consequently $h_nv_n$ is the
difference of two functions which both lie in $\U$ (notice that $f$
lies in $\U$ by Assumption~\ref{basicassumptions}(c)). For the
verification of (\ref{defeqforHD}), we now proceed in two steps.

\textit{Step 1.} To justify the analysis in Step 2 below, we first of all
show that the three integrals
\[
\int|g_{1,f}|(x_1)|\mathrm{d}v_n|(x_1), \qquad \int|g_{2,f}|(x_2)|\mathrm{d}v_n|(x_2), \qquad \int
\int|g(x_1,x_2)| |\mathrm{d}v_n|(x_1)|\mathrm{d}v_n|(x_2)
\]
are finite for all $n\in\N$. For the finiteness of these integrals,
it suffices to show that for every $n \in\N$
\begin{equation}\label{existenceintproofquasi1a}
\int\int|g(x_1,x_2)|\,\mathrm{d}f_n(x_1)\,\mathrm{d}f(x_2)<\infty \quad \mbox{and} \quad \int\int
|g(x_1,x_2)|\,\mathrm{d}f(x_1)\,\mathrm{d}f_n(x_2)<\infty,
\end{equation}
since $|g_{1,f}|\le\int|g(\cdot ,x_2)|\,\mathrm{d}f(x_2)$ and $|g_{2,f}|\le
\int|g(x_1, \cdot)|\,\mathrm{d}f(x_1)$, since $h_n|\mathrm{d}v_n|=\mathrm{d}f_n+\mathrm{d}f$, and since
$f,f_n\in\U$ implies
\[
\int\int|g(x_1,x_2)|\,\mathrm{d}f(x_1)\,\mathrm{d}f(x_2)<\infty
 \quad \mbox{and} \quad   \int\int|g(x_1,x_2)|\,\mathrm{d}f_n(x_1)\,\mathrm{d}f_n(x_2)<\infty.
 \]
  (Notice that
(\ref{existenceintproofquasi1a}) by itself is also needed in Step 2
below.)  We clearly have
\[
\int\int|g(x_1,x_2)|\,\mathrm{d}f(x_1)\,\mathrm{d}f_n(x_2)  \leq  \|\overline{g_{2,f}}\|
_{-\lambda'}\int\phi_{\lambda'}(x_2)\,\mathrm{d}f_n(x_2).
\]
From the second part of Assumption~\ref{basicassumptions}(b) we have
$\|\overline{g_{2,f}}\|_{-\lambda'} < \infty$, and $\int\phi
_{\lambda'}(x_2)\,\mathrm{d}f_n(x_2)<\infty$ holds since $f_n\in\U$. That is, $\|
\overline{g_{2,f}}\|_{-\lambda'}\int\phi_{\lambda
'}(x_2)\,\mathrm{d}f_n(x_2)<\infty$. Similar arguments show that the first
inequality in (\ref{existenceintproofquasi1a}) holds.

\textit{Step 2.} By Step 1 and the triangular inequality we have
\begin{eqnarray}\label
{proofofHDofU}
 &&\biggl|\dot U_f(v)-\frac{U(f+h_nv_n)-U(f)}{h_n} \biggr|\nonumber\\
&& \quad  =   \biggl|-\int v(x_1)\,\mathrm{d}g_{1,f}(x_1)-\int v(x_2)\,\mathrm{d}g_{2,f}(x_2)\nonumber
\\
& & \qquad\hphantom{\biggl|} {} -\frac{1}{h_n} \biggl(\int\int
g(x_1,x_2)\,\mathrm{d}(f+h_nv_n)(x_1)\,\mathrm{d}(f+h_nv_n)(x_2)\nonumber\\
& & \qquad\hphantom{\biggl| {} -\frac{1}{h_n} \biggl(} {} -\int\int g(x_1,x_2)\,\mathrm{d}f(x_1)\,\mathrm{d}f(x_2) \biggr) \biggr|\\
&&  \quad   \le \sum_{i=1}^2 \biggl|-\int v(x_i)\,\mathrm{d}g_{i,f}(x_i)-\int
g_{i,f}(x_i)\,\mathrm{d}v_n(x_i) \biggr| \nonumber\\
& & \qquad {} +   \biggl|h_n\int\int g(x_1,x_2)\,\mathrm{d}v_n(x_1)\,\mathrm{d}v_n(x_2) \biggr|\nonumber
\\
&& \quad  =:  \sum_{i=1}^2S_{1,i}(n)+S_2(n).\nonumber
\end{eqnarray}
In order to show that $S_{1,1}(n)$ converges to zero, we will apply the
integration-by-parts formula~(\ref{defofintegraludv-2}) to
$\int g_{1,f}(x_1)\,\mathrm{d}v_n(x_1)$. At first, we have to make clear that
formula (\ref{defofintegraludv-2}) can be applied, that is, that
the assumptions of Lemma~\ref{intbypart} are fulfilled.

It follows from Step 1 that the second condition in (\ref
{conditionsonintegraludv}) holds true (where $g_{1,f}$ and $v_n$ play
the roles of $u$ and $v$, resp.). Moreover, by the continuity of
$\phi_{-\lambda}$ we have
\begin{eqnarray*}
&&\int|v_n(x_1-)| |\mathrm{d}g_{1,f}|(x_1)\\
&& \quad  =  \int|v_n(x_1-)\phi_{\lambda}(x_1-)\phi_{-\lambda}(x_1-)|
|\mathrm{d}g_{1,f}|(x_1)\\
&& \quad  =  \int|v_n(x_1-)\phi_{\lambda}(x_1-)|\phi_{-\lambda}(x_1)
|\mathrm{d}g_{1,f}|(x_1)\\
&& \quad  \le \|v_n\|_{\lambda}\int\phi_{-\lambda}(x_1) |\mathrm{d}g_{1,f}|(x_1).
\end{eqnarray*}
By Assumption~\ref{basicassumptions}(b) and the fact that $v_n\in\D
_{\lambda}$, the latter bound is finite, so that also the first
condition in (\ref{conditionsonintegraludv}) holds true. We finally
note that $\lim_{|x_1|\to\infty} v_n(x_1)g_{1,f}(x_1)=0$. Indeed: On
one hand, $|g_{1,f}(x_1)\phi_{-\lambda'}(x_1)|$ is bounded above
uniformly in $x_1$ by Assumption~\ref{basicassumptions}(b) and
Remark~\ref{basicassumptionsremark}(b)$'$. On the other hand,
$|v_n(x_1)\phi_{\lambda'}(x_1)|$ converges to $0$ as $|x_1|\to\infty
$ because $|v_n(x_1)\phi_{\lambda}(x_1)|$ is bounded above uniformly
in $x_1$ (recall $\lambda> \lambda'$). That is, the assumptions of
Lemma~\ref{intbypart} are indeed fulfilled.

Now, we may apply the integration-by-parts formula (\ref
{defofintegraludv-2}) to obtain
\begin{eqnarray*}
S_{1,1}(n)
& = &  \biggl|-\int v(x_1)\,\mathrm{d}g_{1,f}(x_1)+\int v_n(x_1-)\,\mathrm{d}g_{1,f}(x_1)
\biggr|\\
& \le&  \biggl|\int(v_n-v)(x_1)\,\mathrm{d}g_{1,f}(x_1) \biggr|+ \biggl|\int
(v_n(x_1-)-v_n(x_1) )\,\mathrm{d}g_{1,f}(x_1) \biggr|\\
& \le&  (\|v_n-v\|_\lambda+\|v_n-v\|_\lambda+\|v-v_n\|_\lambda
 )\int\phi_{-\lambda}(x_1)|\mathrm{d}g_{1,f}|(x_1).
\end{eqnarray*}
The latter bound converges to zero by Assumption \ref
{basicassumptions}(b) and $\|v-v_n\|_\lambda\to0$. That is,
\mbox{$S_{1,1}(n)\to0$}. In the same way we obtain $S_{1,2}(n)\to0$.

Thus, it remains to show $S_2(n)\to0$. We will apply the
integration-by-parts formula~(\ref{defofintegraludv-2}) to the
inner integral in $S_2(n)$. So at first we will verify that formula
(\ref{defofintegraludv-2}) can be used, that is, that the
assumptions of Lemma~\ref{intbypart} are fulfilled. By Assumption
\ref{basicassumptions}(a), we have $g_{x_2}\in\BV_{\loc,\dis}$,
and as mentioned above we also have $v_n\in\BV_{\loc,\dis}$. Further,
the integrals $\int g(x_1,x_2)\,\mathrm{d}f(x_1)$ and $\int
g(x_1,x_2)\,\mathrm{d}f_n(x_1)$ exist by the fact that $f_n, f \in\U$ and
Fubini's theorem. This and the representation $v_n=(f_n-f)/h_n$
imply $\int|g_{x_2}(x_1)| |\mathrm{d}v_n|(x_1)<\infty$, that is, that the
second condition in (\ref{conditionsonintegraludv}) holds true.
Moreover, by the continuity of $\phi_{-\lambda}$ we have as above
\begin{eqnarray*}
\int|v_n(x_1-)| |\mathrm{d}g_{x_2}|(x_1)
& = & \int|v_n(x_1-)\phi_{\lambda}(x_1-)\phi_{-\lambda}(x_1-)|
|\mathrm{d}g_{x_2}|(x_1)\\
& = & \int|v_n(x_1-)\phi_{\lambda}(x_1-)|\phi_{-\lambda}(x_1)|
|\mathrm{d}g_{x_2}|(x_1)\\
& \le& \|v_n\|_{\lambda}\int\phi_{-\lambda}(x_1) |\mathrm{d}g_{x_2}|(x_1).
\end{eqnarray*}
By Assumption~\ref{basicassumptions}(a) and the fact that $v_n
\in\D_{\lambda}$, this bound is finite, so that also the first
condition in (\ref{conditionsonintegraludv}) holds true. We
finally note that $\lim_{|x_1|\to\infty}v_n(x_1)g_{x_2}(x_1)=0$. Indeed:
On one hand, $|g_{x_2}(x_1)\phi_{-\lambda'}(x_1)|$ is bounded above
uniformly in $x_1$ by Assumption~\ref{basicassumptions}(a). On
the other hand, $|v_n(x_1)\phi_{\lambda'}(x_1)|$ converges to $0$ as
$|x_1|\to\infty$ since $|v_n(x_1)\phi_{\lambda}(x_1)|$ is bounded
above uniformly in $x_1$ (recall $\lambda> \lambda'$). That is, the
assumptions of Lemma~\ref{intbypart} are indeed fulfilled.

Now, we may apply the integration-by-parts formula (\ref
{defofintegraludv-2}) to the inner integral in $S_2(n)$ to obtain
\begin{eqnarray*}
S_2(n)
& = &  \biggl|-\int\int v_n(x_1-)\,\mathrm{d}g_{x_2}(x_1)\,\mathrm{d}(f_n-f)(x_2) \biggr|\\
& \le&  \biggl|-\int\int(v_n(x_1-)-v(x_1-))\,\mathrm{d}g_{x_2}(x_1)\,\mathrm{d}(f_n-f)(x_2) \biggr|\\
& &{} +   \biggl|\int\int v(x_1-)\,\mathrm{d}g_{x_2}(x_1)\,\mathrm{d}(f_n-f)(x_2) \biggr|.
\end{eqnarray*}
Since $f_n$ and $f$ generate positive (probability) measures, and $v$
and $\phi_{-\lambda'}$ are continuous, we may continue with
\begin{eqnarray*}
& \le& \|v_n-v\|_{\lambda}\int \biggl(\int\phi_{-\lambda
}(x_1)|\mathrm{d}g_{x_2}|(x_1)
\phi_{-\lambda'}(x_2) \biggr) \phi_{\lambda'}(x_2)\,\mathrm{d}f_n(x_2)\\
& &{} + \|v_n-v\|_{\lambda}\int \biggl(\int\phi_{-\lambda
}(x_1)|\mathrm{d}g_{x_2}|(x_1)
\phi_{-\lambda'}(x_2) \biggr) \phi_{\lambda'}(x_2)\,\mathrm{d}f(x_2)\\
& &{} +   \biggl|\int \biggl(\int v(x_1)\,\mathrm{d}g_{x_2}(x_1) \biggr)\,\mathrm{d}f_n(x_2)-\int
 \biggl(\int v(x_1)\,\mathrm{d}g_{x_2}(x_1) \biggr)\,\mathrm{d}f(x_2) \biggr|\\
& \le& \|v_n-v\|_{\lambda}\int C \phi_{\lambda'}(x_2)\,\mathrm{d}f_n(x_2)
+  \|v_n-v\|_{\lambda}\int C \phi_{\lambda'}(x_2)\,\mathrm{d}f(x_2)\\
& & {}+   \biggl|\int \biggl(\int v(x_1)\,\mathrm{d}g_{x_2}(x_1) \biggr)\,\mathrm{d}f_n(x_2)-\int
\biggl (\int v(x_1)\,\mathrm{d}g_{x_2}(x_1) \biggr)\,\mathrm{d}f(x_2) \biggr|\\
& =: & S_{2,1}(n)+S_{2,2}(n)+S_{2,3}(n)
\end{eqnarray*}
with $C:=\sup_{x_2}\int\phi_{-\lambda}(x_1)|\mathrm{d}g_{x_2}|(x_1) \phi
_{-\lambda'}(x_2)$ (which is finite by the second part of Assumption~\ref{basicassumptions}(a)). By Lemma~\ref{convofintegrals}, which
can be
applied due to Assumption~\ref{basicassumptions}(a), and the
facts that $v \in\D_{\lambda}$, $\|f_n-f\|_{\lambda} \rightarrow
0$, and that $\int
\phi_{\lambda'}(x_2)\,\mathrm{d}f(x_2)$ and $\int
\phi_{\lambda'}(x_2)\,\mathrm{d}f_n(x_2)$ exist, the summand $S_{2,3}(n)$
converges to $0$. Since $\|v_n-v\|_{\lambda}\to0$, and since
$\int\phi_{\lambda'}(x_2)\,\mathrm{d}f(x_2)$ is finite because $f \in\U$, we
also obtain $S_{2,2}(n)\to0$. It remains to show $S_{2,1}(n)\to0$.
As $\|v_n-v\|_{\lambda}\to0$, it suffices to show that
$\int\phi_{\lambda'}(x_2)\,\mathrm{d}f_n(x_2)$ is uniformly bounded from above.
The latter follows from the finiteness of
$\int\phi_{\lambda'}(x_2)\,\mathrm{d}f(x_2)$ and Lemma~\ref{convofintegrals}
which is applicable since we clearly have
$\phi_{\lambda'} \in\D_{-\lambda'}$, and for every $n\in\N$ the
integral $\int\phi_{\lambda'}(x_2)\,\mathrm{d}f_n(x_2)$ exists due to $f_n
\in\U$. This proves the claim of Lemma~\ref{HDofU}.
\end{pf}

\begin{remarknorm}\label{remarkond3}
We note that the proof of Lemma~\ref{HDofU} basically applies also to
V-functionals of the shape
$U(F)=\int\cdots\int g(x_1,\ldots,x_d)\,\mathrm{d}F(x_1)\cdots \mathrm{d}F(x_d)$
with arbitrary \mbox{$d\ge2$}, provided
Assumptions~\ref{basicassumptions}(a)--(b) (which ensure the
quasi-Hadamard differentiability of $U$ in the case $d=2$) are modified
suitably and the definition of $\dot U_f$ in (\ref
{DefHadamard-functional}) is replaced by $\dot U_f(v):=-\sum_{i=1}^d
\int v(x)\,\mathrm{d}g_{i,f}(x)$ with $g_{i,f}(x_i):=\int\cdots\int
g(x_1,\ldots,x_d)\,\mathrm{d}f(x_1)\cdots \allowbreak\mathrm{d}f(x_{i-1})\,\mathrm{d}f(x_{i+1})\cdots \mathrm{d}f(x_d)$.
In particular, Theorem~\ref{maintheorem} then still holds for such
general V-functionals. Let us exemplify the validity of the analogue of
Lemma~\ref{HDofU} for the case $d=3$. To do so, we let $\M_{(\lambda
,\lambda)}$ be the space of all measurable functions $h\dvtx  \R^2 \to\R
$ such that $\sup_{x_1,x_2}|h(x_1,x_2) \phi_{\lambda}(x_1)\phi
_{\lambda}(x_2)|$ is finite. To ensure the existence of the integrals
as in Step 1 in the above proof, it is sufficient to require that the
functions $\overline{g_{i,j,f}}(x_i,x_j):=\int
|g(x_1,x_2,x_3)|\,\mathrm{d}f(x_k)$, $i,j,k \in\{1,2,3\}$, $i < j$, $k \neq i, k
\neq j$, are in $\M_{(-\lambda',-\lambda')}$, and that the functions
$\overline{g_{i,f}}(x_i):=\int|g(x_1,x_2,x_3)|\,\mathrm{d}f(x_j)\,\mathrm{d}f(x_{k})$,
$i,j,k \in\{1,2,3\}$ pairwise disjoint, lie in $\D_{-\lambda'}$
(cf. the second part of Assumption~\ref{basicassumptions}(b)). Then
Step 1 still holds. Let us turn to Step 2 in the above proof. In (\ref
{proofofHDofU}), we now obtain the bound
\begin{eqnarray*}
S_1(n)+S_2(n)+S_3(n)
& := & \sum_{i=1}^3 \biggl|-\int v(x_i)\,\mathrm{d}g_{i,f}(x_i)-\int
g_{i,f}(x_i)\,\mathrm{d}v_n(x_i) \biggr|\\
& &{} +  h_n\sum_{i,j=1:i<j}^3 \biggl|\int\int
g_{i,j,f}(x_i,x_j)\,\mathrm{d}v_n(x_i)\,\mathrm{d}v_n(x_j) \biggr|\\
& & {}+  h_n^2  \biggl|\int\int\int
g(x_1,x_2,x_3)\,\mathrm{d}v_n(x_1)\,\mathrm{d}v_n(x_2)\,\mathrm{d}v_n(x_3)
\biggr|,
\end{eqnarray*}
where $g_{i,j,f}(x_i,x_j):=\int g(x_1,x_2,x_3)\,\mathrm{d}f(x_k)$, $i,j,k \in\{
1,2,3\}$, $i < j$, $k \neq i, k \neq j$. To obtain $S_1(n)\to0$, it
suffices to assume that the functions $g_{i,f}$ satisfy the first part
of Assumption~\ref{basicassumptions}(b). To ensure that
$h_n^{-1}S_2(n)$ is bounded above, it suffices to assume that, similar
to the case $d=2$, the functions $g_{i,j,f}$ satisfy Assumption \ref
{basicassumptions}(a) (with $g$ replaced by $g_{i,j,f}$). Assuming
that for every fixed $x_2,x_3$ the function $g_{x_2,x_3}(\cdot
):=g(\cdot,x_2,x_3)$, lies in $\BV_{\loc} \cap\D_{-\lambda'}$,
and that $(x_2,x_3) \mapsto\int\phi_{-\lambda}(x_1)
|\mathrm{d}g_{x_2,x_3}|(x_1)$ lies in $\M_{(-\lambda',-\lambda')}$ (cf.
Assumption~\ref{basicassumptions}(a)), ensures that $h_n^{-2}S_3(n)$
is bounded above. Thus, $S_1(n)+S_2(n)+S_3(n)\to0$.

Finally, we note that the case $d=1$ is even easier. Here, we only need
to assume $g\in\BV_\loc\cap\D_{-\lambda'}$ (instead of
Assumptions~\ref{basicassumptions}(a)--(b)) and to replace (\ref
{DefHadamard-functional}) by $\dot U_f(v):=-\int v(x)\,\mathrm{d}g(x)$.
\end{remarknorm}


\begin{appendix}
\renewcommand{\theequation}{\arabic{equation}}
\section{Jordan decomposition of functions in $\BV_{\loc}$}\label
{Jordandecomposition}

\setcounter{equation}{17}
Recall that for $\psi\in\BV_{\loc}$ and $c\in\R$, the Jordan
decomposition of $\psi$ centered at $c$,
\begin{eqnarray}\label{eqJordanwitha}
\psi=\psi(c)+\psi_c^+-\psi_c^-,
\end{eqnarray}
is characterized as follows: $\psi_c^+$ and $\psi_c^-$ are the unique
nondecreasing functions satisfying
\begin{eqnarray}
\psi_c^+(x)&=&V^+([c,x],\psi), \qquad  \psi_c^-(x)=V^-([c,x],\psi)
 \qquad \forall x\ge c,\label{jorddecomp0+}\\
\psi_c^+(x)&=&-V^+([x,c],\psi), \qquad  \psi_c^-(x)=-V^-([x,c],\psi
)  \qquad \forall x<c,\label{jorddecomp0-}
\end{eqnarray}
where $V^+([a,b],\psi)$ and $V^-([a,b],\psi)$ denote the positive and
the negative variation of~$\psi$ on the interval $[a,b]$,
respectively. For details see, for example, \cite{Kallenberg2002}, page 34. In our applications, we are mainly concerned with
the positive measures $\mathrm{d}\psi_c^+$ and $\mathrm{d}\psi_c^-$ induced by $\psi
_c^+$ and $\psi_c^-$, respectively (provided  $\psi_c^+$ and  $\psi_c^-$ are right-continuous). The following lemma shows that
$\mathrm{d}\psi_c^+$ and $\mathrm{d}\psi_c^-$ are independent of $c$, although $\psi
_c^+$ and $\psi_c^-$ typically do depend on $c$. In particular, the
definition $|\mathrm{d}\psi|:=\mathrm{d}\psi_c^++\mathrm{d}\psi_c^-$ of the absolute value
measure $|\mathrm{d}\psi|$ is independent of $c$.

\begin{lemma}\label{independenceofc}
Let $\psi\in\BV_{\loc}$ and $c\in\R$. Then $\psi_c^+$, $\psi
_c^-$ differ from $\psi_0^+$, $\psi_0^-$ only by constants
$K_c^+,K_c^-$, respectively. In particular, the positive measures
$\mathrm{d}\psi_c^+$ and $\mathrm{d}\psi_c^-$ are independent of $c$.
\end{lemma}

\begin{pf}
Let $c>0$. Then, in view of (\ref{jorddecomp0+})--(\ref
{jorddecomp0-}), we have
\[
\psi_0^+(x)=V^+([0,x],\psi)=V^+([0,c],\psi)+V^+([c,x],\psi
)=V^+([0,c],\psi)+\psi_c^+(x)
\]
for $x\in(c,\infty)$, and similar we obtain $\psi
_0^+(x)=V^+([0,c],\psi) + \psi_c^+(x)$ for the cases $x\in[0,c]$ and
$x\in(-\infty,0)$. That is, $\psi_c^+=\psi_0^++K_c^+$ for some
constant $K_c^+$. Analogously, we obtain $\psi_c^+=\psi_0^++K_c^+$
for $c\le0$, and $\psi_c^-=\psi_0^-+K_c^-$ for $c\le0$ as well as $c>0$.
\end{pf}


\section{Integration theoretical auxiliaries}\label{appendixintegration}

\setcounter{equation}{20}
Recall our convention $\int=\int_{(-\infty,\infty)}$ and that $\BV
_{\loc,\dis}$ denotes the space of all  c\`adl\`ag  functions in $\BV
_{\loc}$. 

\begin{lemma}\label{intbypart}
Let $u,v\in\BV_{\loc,\dis}$ such that $\lim_{x\to\pm\infty
}u(x)v(x)=c_\pm$ for some constants $c_+, c_-\in\R$. Then, if
\begin{equation}\label{conditionsonintegraludv}
\int|v(x-)| |\mathrm{d}u|(x)<\infty  \quad \mbox{and} \quad  \int|u(x)|
|\mathrm{d}v|(x)<\infty,
\end{equation}
we have the integration-by-parts formula
\begin{equation}\label{defofintegraludv-2}
\int u(x)\,\mathrm{d}v(x) = c_+-c_--\int v(x-)\,\mathrm{d}u(x).\vspace*{-3pt}
\end{equation}
\end{lemma}

\begin{pf}
If $-\infty<a<b<\infty$, then one can proceed as in the proof of
Theorem II.6.11 in~\cite{Shiryaev1996} to obtain
\begin{equation}\label{defofintegraludv-2-proof}
\int_{(a,b]}u(x)\,\mathrm{d}v(x) = u(b)v(b)-u(a)v(a)-\int_{(a,b]}v(x-)\,\mathrm{d}u(x),
\end{equation}
because $\int_{(a,b]}|v(x-)| |\mathrm{d}u|(x)<\infty$ and $\int
_{(a,b]}|u(x)| |\mathrm{d}v|(x)<\infty$. 
Now, choosing sequences $(a_n),(b_n)\subset(-\infty,\infty)$ with
$a_n\downarrow-\infty$ and $b_n\uparrow\infty$, the statement of
the lemma follows from~(\ref{defofintegraludv-2-proof}), the
continuity from below of the finite measures $\int_. u^+(x)\,\mathrm{d}v^+(x)$,
$\int_. u^-(x)\,\mathrm{d}v^+(x), \ldots$ on $(-\infty,\infty)$, and the
assumption $\lim_{x\to\pm\infty}u(x)v(x)=c_\pm$.\vspace*{-3pt}
\end{pf}

Next, we give a sort of Helly--Bray theorem. Recall that $\BV_{1,\dis
}$ denotes the space of all  c\`adl\`ag  functions on $\oR$ with
variation bounded by $1$.\vspace*{-3pt}

\begin{lemma}\label{convofintegrals}
Let $\lambda>\lambda'\ge0$, let $\psi\in\D_{-\lambda'}$ and suppose
that $f,f_1,f_2,\ldots\in\BV_{1,\dis}$ are nondecreasing and satisfy
$\lim_{n\to\infty}\|f_n-f\|_\lambda=0$. Let
$\int\phi_{\lambda'}(x)\,\mathrm{d}f(x) < \infty$ and
$\int\phi_{\lambda'}(x)\,\mathrm{d}f_n(x) < \infty$ for every $n \in\N$. Then
the integrals $\int\psi(x)\,\mathrm{d}f(x)$ and $\int\psi(x)\,\mathrm{d}f_n(x)$ exist and
we have
\[
\lim_{n\to\infty}\int\psi(x)\,\mathrm{d}f_n(x)=\int\psi(x)\,\mathrm{d}f(x).\vspace*{-3pt}
\]
\end{lemma}

\begin{pf}
The first claim follows from
\begin{eqnarray*}
\int|\psi(x)|\,\mathrm{d}f(x)=\int|\psi(x) \phi_{\lambda'}(x) \phi
_{-\lambda'}(x)|\,\mathrm{d}f(x) \leq\|\psi\|_{-\lambda'} \int\phi_{\lambda
'}(x)\,\mathrm{d}f(x)
\end{eqnarray*}
and the analogous bound for $\int|\psi(x)|\,\mathrm{d}f_n(x)$, $n \in\N$.

Now let us turn to the second claim. Since $\psi\phi_{-\lambda'}$ is
a bounded  c\`adl\`ag  function on the \textit{compact} interval $\oR$, we
may and do choose for each $\varepsilon> 0$ a step function
$\widetilde\psi_\varepsilon\in\D$ with a~finite number of jumps
and satisfying $\|\psi\phi_{-\lambda'}-\widetilde\psi_\varepsilon
\|_\infty\leq\varepsilon$. For $\psi_\varepsilon:= \widetilde\psi
_\varepsilon\phi_{\lambda'}$, we thus have $\|\psi-\psi
_\varepsilon\|_{-\lambda'}\le\varepsilon$. Of course,
\begin{eqnarray}\label{convofintegrals-proof-eq-1}
 &&\biggl| \int\psi(x)\,\mathrm{d}f_n(x)-\int\psi(x)\,\mathrm{d}f(x)
\biggr|\nonumber\\[-2pt]
&& \quad  \le  \biggl| \int\psi(x)\,\mathrm{d}(f_n-f)(x) - \int\psi_\varepsilon(x)\,\mathrm{
d}(f_n-f)(x) \biggr|\nonumber
\\[-9pt]
\\[-9pt]
& & \qquad {} + \biggl|\int\psi_\varepsilon(x)\,\mathrm{ d}(f_n-f)(x) \biggr|\nonumber\\
&& \quad  =:  S_1(n,\varepsilon)+S_2(n,\varepsilon).
\nonumber
\end{eqnarray}
For the first summand, we obtain
\begin{eqnarray}\label{eqginftyreplacedbystep}
S_1(n,\varepsilon)
& = &  \biggl| \int\phi_{-\lambda'}(x) \phi_{\lambda'}(x) \psi
(x)\,\mathrm{d}(f_n-f)(x)\nonumber\\
& &\hphantom{\biggl|}{} - \int\phi_{-\lambda'}(x) \phi_{\lambda'}(x) \psi_\varepsilon
(x)\,\mathrm{ d}(f_n-f)(x)  \biggr| \nonumber\\
& \le& \biggl (\int\phi_{\lambda'}(x)\,\mathrm{d}f_n(x) + \int\phi_{\lambda
'}(x)\,\mathrm{d}f(x) \biggr)\|\psi-\psi_\varepsilon\|_{-\lambda'} \\
& \le&  \biggl(\int\phi_{\lambda'}(x)\,\mathrm{d}f_n(x) + \int\phi_{\lambda
'}(x)\,\mathrm{d}f(x)  \biggr) \varepsilon\nonumber\\
& \le& C \varepsilon\nonumber
\end{eqnarray}
for some finite constant $C>0$ being independent of $n$ and
$\varepsilon$. For the last step, we used the assumption $\int
\phi_{\lambda'}(x)\,\mathrm{d}f(x) < \infty$ and the fact that
$\sup_{n\in\N}\int\phi_{\lambda'}(x)\,\mathrm{d}f_n(x) < \infty$. The latter
fact is not completely obvious, so that we give the details: Because
of $\int\phi_{\lambda'}(x)\,\mathrm{d}f(x) < \infty$, it is clearly
sufficient to show that $\sup_{n\in\N}|\int\phi_{\lambda'}(x)\,\mathrm{
d}(f-f_n)(x)|$ is bounded above by some finite constant. By our
assumptions and the bound (\ref{eqfnintegralfinite}) below, we
can apply the integration by parts formula (\ref{defofintegraludv-2})
to the functions $f-f_n$ and $\phi_{\lambda'}$ to obtain
\[
 \biggl| \int\phi_{\lambda'}(x)\,\mathrm{ d}(f-f_n) \biggr| \le 2\|f-f_n\|
_{\lambda'} +  \biggl| \int(f-f_n)(x-)\,\mathrm{d}\phi_{\lambda'}(x)  \biggr|.
\]
By our assumptions, the first summand tends to $0$ since $\|f_n-f\|
_{\lambda'}\leq\|f_n-f\|_{\lambda}$. The second summand is less than
or equal to $\int|(f-f_n)(x-)|  |\mathrm{d}\phi_{\lambda'}|(x)$ and we have
\begin{eqnarray}\label{eqfnintegralfinite}
\int|(f-f_n)(x-)| |\mathrm{d}\phi_{\lambda'}|(x)
& = & \int |\phi_{\lambda}(x) (f-f_n)(x-)| \phi_{-\lambda}(x)
|\mathrm{d}\phi_{\lambda'}|(x) \nonumber
\\[-9pt]
\\[-9pt]
& \leq& 2 \|f-f_n\|_\lambda  \int_{\mathbb{R}_+} \phi_{-\lambda
}(x)\,\mathrm{d}\phi_{\lambda'}(x).
\nonumber
\end{eqnarray}
Since $\|f-f_n\|_\lambda\rightarrow0$ by assumption, and $\int
_{\mathbb{R}_+} \phi_{-\lambda}(x)\,\mathrm{d}\phi_{\lambda'}(x)<\infty$ by
$\lambda>\lambda'\ge0$, the left-hand side of (\ref
{eqfnintegralfinite}) converges to $0$. In particular, the left-hand
side of (\ref{eqfnintegralfinite}) is bounded above uniformly in $n$.
This completes the proof of (\ref{eqginftyreplacedbystep}).

Now, the second claim of the lemma would follow from (\ref
{convofintegrals-proof-eq-1}) and (\ref{eqginftyreplacedbystep}) if we
could show that $S_2(n,\varepsilon)$ converges to $0$ as $n\to\infty
$ uniformly in $\varepsilon\in(0,1]$. By our assumptions and formula
(\ref{eqlaststephellybray}) below, we can apply the integration by
parts formula (\ref{defofintegraludv-2}) to obtain
\begin{eqnarray*}
S_2(n,\varepsilon)
& = &  \biggl| \int\psi_\varepsilon(x) \phi_{\lambda'}(x) \phi
_{-\lambda'}(x) \,\mathrm{d}(f_n-f)(x)  \biggr| \\
&\le& 2\|\psi_\varepsilon\|_{-\lambda'} \|f_n-f\|_{\lambda'} +
 \biggl| \int(f_n-f)(x-)\,\mathrm{d}\psi_\varepsilon(x)  \biggr|\\
&\le& 2(\|\psi_\varepsilon-\psi\|_{-\lambda'}+\|\psi\|_{-\lambda
'}) \|f_n-f\|_{\lambda'} +  \biggl| \int(f_n-f)(x-)\,\mathrm{d}\psi_\varepsilon
(x)  \biggr|.
\end{eqnarray*}
The first summand converges to $0$ by our assumptions and $\|\psi
_\varepsilon-\psi\|_{-\lambda'}\le\varepsilon\le1$. Furthermore,
the second summand is less than or equal to $\int|(f_n-f)(x-)|
|\mathrm{d}\psi_\varepsilon|(x)$. Recalling $\psi_\varepsilon=\widetilde
\psi_\varepsilon\phi_{\lambda'}$ and that $\widetilde\psi
_\varepsilon$ is a step function with a finite number of jumps, we now obtain
\begin{eqnarray}\label
{eqlaststephellybray}
  &&\int|(f_n-f)(x-)|  |\mathrm{d}\psi_\varepsilon|(x) \nonumber\\
&& \quad  \le \|\widetilde\psi_\varepsilon\|_\infty \int|(f_n-f)(x-)|
|\mathrm{d}\phi_{\lambda'}|(x)\nonumber
\\[-8pt]
\\[-8pt]
&& \quad  =  \|\widetilde\psi_\varepsilon\|_\infty \int|(f_n-f)(x-)\phi
_{\lambda}(x)|  \phi_{-\lambda}(x) |\mathrm{d}\phi_{\lambda'}|(x)
\nonumber\\
&& \quad  \leq 2(\|\psi\phi_{-\lambda'}\|_\infty+1) \|f_n-f\|_\lambda
\int_{\mathbb{R}_+} \phi_{-\lambda}(x)\,\mathrm{d}\phi_{\lambda'}(x),
\nonumber
\end{eqnarray}
and this expression converges to $0$ because $\|f_n-f\|_{\lambda}\to
0$ and $\lambda>\lambda'\ge0$. That is, $S_2(n,\varepsilon)$ indeed
converges to $0$ as $n\to\infty$ uniformly in $\varepsilon\in(0,1]$.
\end{pf}
\end{appendix}

\section*{Acknowledgements}
The authors wish to thank an Associate Editor and the reviewers for
their very careful reading and for useful hints and comments.

%

\printhistory

\end{document}